\documentclass[12pt]{article}
\usepackage{amsxtra,amssymb,amsthm,amsmath,latexsym}

\theoremstyle{plain}
\newtheorem{theorem}{Theorem}

\newtheorem{lemma}{Lemma}

\newcommand{\refT}[1]{Theorem~\ref{T:#1}}
\newcommand{\refS}[1]{Section~\ref{S:#1}}

\newcommand{\refL}[1]{Lemma~\ref{L:#1}}

\def\qed{{\hfill $\Box$}}
\def\ve{\varepsilon}

\def\R{{\mathbb R}}

\def\calL{{\mathcal L}}
\def\calM{{\mathcal M}}
\def\calN{{\mathcal N}}

\def\oH{{\overset{\circ}{H}}}
\def\oH1{{\overset{\circ}{H}\kern-.02in{}^1}}
\def\l{\ell}

\def\bee{\begin{equation*}}
\def\eee{\end{equation*}}
\def\be{\begin{equation}}
\def\ee{\end{equation}}

\begin{document}
\title{An inverse problem with data on the part of the boundary  }

\author{A.G. Ramm\\
 Mathematics Department, Kansas State University, \\
 Manhattan, KS 66506-2602, USA\\
ramm@math.ksu.edu }

\date{}
\maketitle\thispagestyle{empty}

\begin{abstract}

Let $u_t=\nabla^2 u-q(x)u:=Lu$ in $D\times [0,\infty)$, where 
$D\subset R^3$ is a bounded domain with a smooth
connected boundary $S$, and $q(x)\in L^2(S)$ is a 
real-valued function with compact support in $D$. Assume that $u(x,0)=0$, 
$u=0$ on
$S_1\subset S$, $u=a(s,t)$ on $S_2=S\setminus S_1$, where 
$a(s,t)=0$ for $t>T$,  $a(s,t)\not\equiv 0$,
$a\in C([0,T];H^{3/2}(S_2))$ is arbitrary.

Given the extra data $u_N|_{S_2}=b(s,t)$,
for each $a\in C([0,T];H^{3/2}(S_2))$, where $N$ is the outer normal to 
$S$, one can find $q(x)$ uniquely. A similar result is obtained for
the heat equation $u_t=\mathcal{L} u:=
\nabla \cdot (a \nabla u)$.

These results are based on  new versions of Property C.
\end{abstract}

\section{Introduction}\label{S:1}
Let $D\subset\R^3$ be a bounded domain with a
 smooth connected boundary $S$, $S=S_1\cup S_2$, $S_1$ is an open subset 
in $S$, and $S_2$ is the complement of $S_1$ in $S$. 
Consider the problem
\be\label{e1}u_t=\nabla^2 u-q(x)u, \qquad (x,t)\in 
D\times[0,\infty), \ee
\be\label{e2}u(x,0)=0, \qquad u|_{S_1}=0, \qquad u|_{S_2}=a(s,t), \ee
where $q(x)\in L^2(D)$ is a real-valued function, 
$a(s,t)\in C([0,\infty),H^{3/2}(S_2))$, the function $a$ is an arbitrary 
real-valued
function in the above set such that $a\not\equiv 0$, $a=0$ for $t>T$, 
$H^\ell$ is the Sobolev space.

For each $a$ in the above set let the extra data be given:
\be\label{e3} u_N|_{S_2}=b(s,t),\ee
where $N$ is the outer normal to the boundary $S$.

{\it Do the data $\{a(s,t),b(s,t)\}$ $\forall a(s,t)\in 
C([0,\infty),H^{3/2}(S_2))$, $a\not\equiv 0$, $a=0$ for $t>T$, determine 
$q(x)$ uniquely?}

Our main result is a positive answer to this question.

\begin{theorem}\label{T:1}
The data $\{a(s,t),b(s,t)\}_{\forall t\geq0,\ \forall s\in S_2}$,  
given for all $a$ with the above properties, determine a compactly
supported in $D$ real-valued function $q(x)\in L^2(D)$ 
uniquely.
\end{theorem}

Actually we prove a slightly stronger result: the data for $0\leq t\leq T+\ve$ determine $q$ uniquely, 
where $\ve>0$ is an arbitrary small number.

Note that the set $S_2$ can be arbitrary small.

\refT{1} is a multidimensional 
generalization of the author's result from \cite{R1}.
Let 
$$L_jv-\lambda v:=\nabla^2 v-q_j(x)v-\lambda v, \quad j=1,2,$$
and  
$$N_j=\{v:L_j v-\lambda v=0\hbox{\ in \ }D, v|_{S_1}=0\}.$$

The main tool in the proof of Theorem 1 is a 
new version of 
property C.  Originally this property was
introduced by the author in 1986  for the products of solutions
to homogeneous linear partial differential equations in the case
when these solutions did not satisfy any boundary conditions (see 
[3]).

\begin{theorem}\label{T:2}
The set $\{v_1,v_2\}$ for all $v_1\in N_1$ and all $v_2\in N_2$ is 
complete (total) in $L^2(D_1)$, where $D_1\subset D$ is a strictly inner 
subdomain of $D$,
i.e., if $p\in L^2(D_1)$, $p=0$ in $D\setminus D_1$, and $\int_{D_1} p v_1 
v_2\,dx=0$ 
\quad $\forall v_1\in N_1, \forall v_2\in N_2$, 
then $p=0$. 
\end{theorem}

In \refS{2} Theorems 1 and 2 are proved. In \refS{3} the results are 
generalized to the boundary-value problems for the  equation
$u_t=\nabla\cdot(a(x)\nabla u)$.

\section{Proofs}\label{S:2}

\begin{proof}[Proof of \refT{1}]
Let \[v=\int^\infty_0e^{-\lambda t}u\,dt \hbox{\quad and 
\quad} 
A(s,\lambda)=\int^\infty_0e^{-\lambda t}a(s,t)\,dt.\] Taking 
the Laplace 
transform of the relations (1)-(2), we obtain:
\be\label{e4}
Lv-\lambda v=0, \quad Lv=\nabla^2 v-q(x)v, \quad v|_{S_1}=0, 
\quad v|_{S_2}=A(s,\lambda), \ee
and \[v_N=B(s,\lambda)=\int^\infty_0 e^{-\lambda t} 
b(s,t)\,dt.\]

Assume that there are $q_1$ and $q_2$, compactly supported in $D$, which 
generate the same data. 
Let \[L_j v_j=\nabla^2 v_j-q_j v_j, \quad j=1,2,\] and 
\[w=v_1-v_2.\] Then
\be\label{e5}L_1w-\lambda w=p v_2,\quad p=q_1-q_2, \ee
and for any \[\psi\in 
N(L_1-\lambda):=\{\psi:(L_1-\lambda)\psi=0\},\quad 
\psi|_{S_1}=0,\]
 one gets
\be\label{e6} \int_D p v_2\psi\,dx=\int_{S_2} (w_N\psi-\psi_N w)ds=0,\ee
because \[w=\psi=0 \hbox{\quad  on\quad } S_1\] and 
\[w=w_N=0 \hbox{\quad  on\quad } S_2\] by our 
assumptions. By 
\refT{2}, relation \eqref{e6} implies $p=0$.

\refT{1} is proved.
\end{proof}

\begin{proof}[Proof of \refT{2}]
It is proved in \cite{R2} that the set $\{v_1 v_2\}$ for all 
\[v_j\in M_j:=\{v_j:(L_j-\lambda)v_j=0\hbox{\ in\ }D\}\] is 
total in $L^2(D)$.
Therefore it is sufficient to prove that $N_j$ is dense in 
$M_j$ in  $L^2(D_1)\hbox{-norm}$, where $D_1$ is a strictly inner 
subdomain of $D$ out of which both $q_1$ and $q_2$ vanish.

Let us take $j=1$. The proof for $j=2$ is the same. Assume 
the contrary. 
Then there is a $\psi\in M_1$ such that
\be\label{e7}0=\int_{D_1}\psi v\,dx \qquad \forall v\in N_1. \ee
Let $G(x,y)$ solve the problem
\be\label{e8}(L_1-\lambda)G=-\delta(x-y)\hbox{\ in\ }\R^3, 
  \qquad \lim_{|x-y|\to\infty}G(x,y)=0, \qquad G|_{S_1}=0.\ee
Since $G(x,y) \in N_1$  $\forall y\in D_1'=\R^3\setminus D_1$, equation 
\eqref{e7} implies
\be\label{e9} 0=\int_{D_1} \psi(x) G(x,y) dx:=h(y), \qquad \forall y\in 
D_1'.\ee  
We have $h(y)\in H^2_{loc}(\R^3\setminus S_1)$, and $h$ solves the 
elliptic equation 
\[ (\nabla^2 -\lambda)h=0 \hbox{ \quad  in  \quad} D_1',\]
because $q_1=q_2=0$ in $D_1'$.
Therefore, by the uniqueness of the solution to the Cauchy problem for 
elliptic equtions,  
one gets from (9)
the following relations:
\be\label{e10} h=h_N=0 \hbox{\ on\ } \partial D_1, \ee
and, because of \eqref{e8}, one gets
\be\label{e11}
(L_1-\lambda)h=-\psi \hbox{\  in   \ }D_1. \ee
Multiply \eqref{e11} by $\overline{\psi}$, integrate over $D_1$, 
use \eqref{e10}, and get
\be\label{e12} -\int_{D_1} |\psi|^2 dx=\int_{\partial D_1} (h_N 
\overline{\psi}-h\overline{\psi_N})ds=0. \ee

Thus $\psi=0$ in $D_1$ and, therefore, $\psi=0$ in $D$, because $\psi$ solves a 
homogeneous 
linear elliptic equation for which the uniqueness of the solution to the Cauchy 
problem holds.  \refT{2} is proved. \end{proof}

\section{Generalizations }\label{S:3}
Consider now the problem
\be\label{e13} u_t=\nabla\cdot(a(x)\nabla u):=\calL u, \qquad (x,t)\in 
D\times 
[0,\infty), \ee
\be\label{e14} u(x,0)=0, \qquad u|_{S_1}=0, \qquad u|_{S_2}=h(s,t), \ee
and the extra data are
\be\label{e15}a u_N|_{S_2}=z(s,t) \quad \forall s\in S_2, 
\forall t>0.\ee

We assume that \[0<a_0\leq a(x)\leq a_1,\] where $a_0$ and 
$a_1$ are constants, $a\in 
H^3(D)$,
and  prove
the following theorem.
 
\begin{theorem}\label{T:3}
Under the above assumptions, the data $\{h(s,t), z(s,t)\}_{\forall s\in S_2, \forall 
t>0}$ determine 
$a(x)$ uniquely.
\end{theorem}

\begin{proof} Taking the Laplace transform of relations (13)-(15), we 
reduce the problem to
\be\label{e16} (\calL-\lambda)v=0\hbox{\ in\ }D, \quad v|_{S_1}=0, \
v|_{S_2}=H(s,\lambda), \quad a v_N|_{S_2} =Z(s,\lambda), \ee 
where, e.g., \[H(s,\lambda)=\int_0^\infty e^{-\lambda 
t}h(s,t)dt.\]
Assuming that
$a_j,j=1,2,$ generate the same data $\{H(s,\lambda), Z(s,\lambda)\}_{\forall s\in 
S_2, \forall \lambda>0}$, one
derives for \[w=v_1-v_2\] the problem \be\label{e17}
(\calL_1-\lambda)w=\nabla\cdot(p\nabla v_2),\quad p:=a_2-a_1,\quad 
w|_{S}=0,\quad a_1v_{1N}|_{S_2}=a_2v_{2N}|_{S_2}. 
\ee 
Multiply \eqref{e17} by an arbitrary
element of $\calN_1$, where
\[\calN_j=\{\varphi:(\calL_j-\lambda)\varphi=0\hbox{\ in\ }D,\quad
\varphi|_{S_1}=0 \}, \quad j=1,2,\] integrate by parts, and get 
\be\label{e18} -\int_D
p\nabla v_2\nabla \varphi\,dx+\int_S p v_{2N}\varphi\, ds
=\int_S(a_1w_N\varphi-a_1\varphi_N w)ds. \ee 
Using boundary conditions
\eqref{e17}, one gets 
\be\label{e19} \int_D p \nabla\varphi \nabla v_2dx=0
\qquad \forall \varphi\in\calN_1\quad \forall v_2\in \calN_2. \ee

To complete the proof of Theorem 3, we use the following new version of Property C.

\begin{lemma}\label{L:1}
The set $\{\nabla\varphi\cdot\nabla v_2\}_{\forall \varphi\in\calN_1,
\forall u_2\in\calN_2}$ is complete in $L^2(D)$ for all sufficiently large
$\lambda>0$.
\end{lemma}

We prove this lemma below, but first let us explain the 
claim made in the Introduction: 

{\it Claim:} The results remain valid if the data are given not 
for all $t\geq 0$ but for $t\in[0,T+\ve]$, where $\ve>0$ is arbitrarily 
small.

This claim follows from the analyticity with respect to time of the 
solution $u(x,t)$ to problems 
\eqref{e1}-\eqref{e2} and \eqref{e13}-\eqref{e14} in a 
neighborhood of the ray $(T,\infty)$ for an arbitrary small $\ve>0$. 
This analyticity holds if $a(s,t)$ and $h(s,t)$ vanish in the region $t>T$.

\begin{proof}[Proof of \refL{1}]
It was proved in \cite[pp.78-80]{R2}, that the set 
$\{\nabla\psi_1\cdot\nabla\psi_2\}_{\forall\psi_j\in \calM_j}$
is complete in $L^2(D)$, where
$$\calM_j=\{\psi:(\calL_j-\lambda)\psi_j=0\},$$ 
$$\calL_j\psi=\nabla\cdot(a_j\nabla\psi)-\lambda\psi,$$
where
$\lambda=const>0$, $0<c\leq a_j(x)\leq C$, $a_j(x)\in H^3(D)$,
$j=1,2,$ and $H^\l(D)$ is the usual Sobolev space. 

\refL{1} will be proved if we 
prove that any $\psi_j\in \calM_j$ 
can be approximated with an arbitrary small error in the norm, generated 
by the bilinear form
\[\int_D(a_j\nabla\psi_j\nabla\varphi+\lambda\psi_j\varphi)dx\] by an 
element 
$v_j\in \calN_j$. The above norm is equivalent to the norm of $H^1(D)$ due
to the assumption $0<c\leq a_j(x)\leq C, \quad j=1,2$.

Assuming that such an approximation is not possible, we 
can find a $\psi_j\in \calM_j$ 
such that
\be\label{e20}
0=\int_D [a_j(x)\nabla\psi_j\nabla G_j(x,y)+\lambda\psi_j G(x,y)]dx 
\qquad \forall y\in D', \ee
because $G_j(x,y)\in \calN_j$ for any $y\in D'$.
Integrating by parts in \eqref{e20} and using the relation $G_j\in\calN_j$, 
one gets
\be\label{e21}
0=\int_S a_j\psi_j G_{jN}(s,y)ds:=h(y), \qquad \forall y\in D'. \ee
Denote \[a_j\psi_j=\varphi,\quad G_j=G,\quad \calL_j=\calL.\]
 Since
\be\label{e22}
(\calL-\lambda)G=-\delta(x-y) \hbox{\  in \ } \R^3\setminus S_1 \hbox{  
and \ } G|_{S_1}=0,\ee
one can derive the relation
\be\label{e23}\lim_{y\to p\in S_1} G_N(s,y)=\frac{\delta(s-p)}{a(s)}, \ee
where $\delta(s-p)$ is the delta-function on $S_1$.
We prove this relation later, but assuming that \eqref{e23} holds we conclude 
from \eqref{e21} that
\be\label{e24} a_j\psi_j:=\varphi=0\hbox{\quad on\quad \ }S_1. \ee
Consequently
\be\label{e25}
h(y)=\int_{S_2} \varphi(s) G_N(s,y)ds=0 \qquad \forall y\in D'.\ee
It follows from \eqref{e25} that
\be\label{e26} (\calL-\lambda)h=0 \hbox{\quad in\quad \ }  \R^3\setminus 
S_2, 
\qquad 
h=0\hbox{\ in\ } D'. \ee
This and  the uniqueness of the solution of the Cauchy problem 
for elliptic equation \eqref{e25},
imply
\be\label{e27} h=0 \hbox{\ in  \ } \R^3\setminus S_2. \ee
From \eqref{e27} and the  jump relation for the double-layer 
potential \eqref{e25}, we conclude that
\be\label{e28} \varphi=0 \hbox{\  on  \ } S_2. \ee
From \eqref{e24} and \eqref{e28} it follows that $\varphi=0$ on $S$. 
Therefore $\psi_j=0$ on $S$ because $a_j>0$. Thus,  $\psi_j\in \calM_j$ 
and $\psi_j=0$ on $S$. This implies  $\psi_j=0$ in $D$ for 
all sufficiently large $\lambda>0$. 

Thus, \refL{1} is proved if \eqref{e23} is established.

Let us prove \eqref{e23}. Denote by $g$ the (unique) solution to the problem
\be\label{e29} (\calL-\lambda)g=-\delta(x-y)\hbox{\ in\ }\R^3, \quad \lim_{|x-y|\to 
\infty}g(x,y)=0. \ee
Using Green's formula, we obtain from \eqref{e22} and \eqref{e29} the 
relation
\be\label{e30}
G(x,y)=g(x,y)-\int_{S_1} g(x,s) a(s) G_N(s,y)ds, 
\quad x,y\in S'_1:=\R^3\setminus S_1. \ee
Taking $y\to p\in S_1$ and using the boundary condition \eqref{e22}, we derive:
\be\label{e31} 
g(x,p)=\lim_{y\to p} \int_{S_1} g(x,s) a(s) G_N(s,y) ds \qquad \forall x\in S'_1.\ee
Since the set $\{g(x,t)\}_{\forall x\in S'_1}$ is dense in $L^2(S_1)$, formula \eqref{e31} 
implies \eqref{e23}. 

Let us finally check the claim that the set $\{g(x,p)\}_{\forall x\in S'_1}$
is dense in $L^2(S_1)$. If it is not dense, then there is an $f\in L^2(S_1)$ such that
\be\label{e32} 0=\int_{S_1} f(p)g(x,p)dp:=W(x), \qquad \forall x\in S'_1. 
\ee
The integral in \eqref{e32} is a simple-layer potential $W(x)$ the density 
$f(p)$ of which
must vanish because of the jump formula for the normal derivatives of $W$ 
across $S_1$. 
Thus the claim is proved.
The proof of \refL{1} is complete.  Therefore  \refT{3} is proved.
\end{proof}

\renewcommand{\qed}{}
\end{proof}

{\bf Remark.} If the conclusion of Lemma 1 remains valid  for the set
$\{\nabla\varphi\cdot\nabla v_2\}_{\forall \varphi\in\calN_1, u_2
\in\calN_2}$, where $u_2$ is a single element of $\calN_2$,
possibly chosen in a special way, then the conclusion of Theorem 1 will be 
established for the data which is a single pair of data
$\{a(s,t),b(s,t)\}_{\forall t\geq0,\ \forall s\in S_2}$,
where $a(s,t)\not\equiv 0$ is some function.

\end{document}